\documentclass[12pt]{amsart}
\usepackage{amssymb,amsmath,mathrsfs}
\usepackage[ps2pdf,colorlinks=true,urlcolor=blue,
citecolor=red,linkcolor=blue,linktocpage,pdfpagelabels,bookmarksnumbered,bookmarksopen]{hyperref}
\usepackage{epsfig,graphicx,color,mathrsfs}
\usepackage[english]{babel}

\usepackage[left=2.8cm,right=2.8cm,top=2.75cm,bottom=2.75cm]{geometry}

\newcommand{\R}{{\mathbb R}}

\newcommand{\RN}{{\mathbb{R}^{N}}}

\newcommand{\tss}{\bar t}

\newcommand{\abs}[1]{\lvert#1\rvert}
\newcommand{\norm}[2][]{\|#2\|_{#1}}

\renewcommand{\H}{{\mathbb H}}

\newcommand{\rmd}{\mathrm{d}}

\newcommand{\iu}{{\rm i}}
\newcommand{\teta }{\theta }
\newcommand{\vep}{\varepsilon}
\newcommand{\eps}{\varepsilon}

\numberwithin{equation}{section}
\newtheorem{theorem}{Theorem}[section]
\newtheorem{proposition}[theorem]{Proposition}
\newtheorem{lemma}[theorem]{Lemma}
\newtheorem{ana}[theorem]{Analytical Property}

\theoremstyle{definition}

\newcommand{\brm}{\begin{remark}\rm}
\newcommand{\erm}{\end{remark}}
\newcommand{\brms}{\begin{remark}\rm}
\newcommand{\erms}{\end{remark}}
\newcommand{\bte}{\begin{theorem}}
\newcommand{\ete}{\end{theorem}}
\newcommand{\bpr}{\begin{proposition}}
\newcommand{\epr}{\end{proposition}}
\newcommand{\ble}{\begin{lemma}}
\newcommand{\ele}{\end{lemma}}
\newcommand{\beq}{\begin{equation}}
\newcommand{\eeq}{\end{equation}}
\newcommand{\bdm}{\begin{displaymath}}
\newcommand{\edm}{\end{displaymath}}
\numberwithin{equation}{section}

\newcommand{\bos}{\begin{remark}\rm}
\newcommand{\eos}{\end{remark}}

\newcommand{\ben}{\begin{enumerate}}
\newcommand{\een}{\end{enumerate}}

\newcommand{\be}{\begin{equation}}
\newcommand{\ee}{\end{equation}}

\title[Numerical computation of soliton dynamics for NLS equations]{Numerical computation of
soliton dynamics \\ for NLS equations in a driving potential}

\author[M.\ Caliari]{Marco Caliari}
\author[M.\ Squassina]{Marco Squassina}

\address{Dipartimento di Informatica
\newline\indent
Universit\`a degli Studi di Verona
\newline\indent
C\'a Vignal 2, Strada Le Grazie 15, I-37134 Verona, Italy}
\email{marco.caliari@univr.it}
\email{marco.squassina@univr.it}

\thanks{The first author was partially supported by the GNCS ``Programma Giovani Ricercatori''.}

\thanks{The second author was partially supported by the
Italian PRIN Project 2007: {\em Metodi Variazionali e Topologici
nello Studio di Fenomeni non Lineari}}

\begin{document}



\begin{abstract}
We provide some numerical computations for the soliton dynamics
of the nonlinear Schr\"odinger equation with an external potential.
After computing the ground state
solution $r$ of a related elliptic equation we show that, in the semi-classical
regime, the center of mass of the solution with initial datum modelled on $r$ is driven by 
the solution of $\ddot x=-\nabla V(x)$. Finally, we 
provide some examples and analyze the numerical errors in the two dimensional case
when $V$ is an harmonic potential.
\end{abstract}
\maketitle

\section{Introduction}

\subsection{Soliton dynamics behaviour}
The goal of this paper is to provide
a numerical investigation of the so-called {\em soliton dynamics} behaviour for the nonlinear Schr\"odinger 
equation with an external time independent smooth potential $V$
\begin{equation}
\label{probMF}
\tag{$P$}
\begin{cases}
\iu\eps\partial_t\phi_\eps=-\frac{\eps^2}{2}\Delta \phi_\eps
+V(x)\phi_\eps-|\phi_\eps|^{2p}\phi_\eps,\quad & \text{$x\in\R^N,\, t>0$}, \\
\noalign{\vskip4pt}
\phi_\eps(x,0)=\phi_0(x),\quad & \text{$x\in\R^N$},
\end{cases}
\end{equation}
that is the qualitative behaviour of the solution $\phi_\eps(t)$ in the 
semi-classical regime, namely for $\eps$ (which of course plays the r\^ole of 
Planck's constant) going to zero, 
by taking as initial datum a (bump like) function of the form
\begin{equation}
\label{initialD}
\tag{$I$}
\phi_0(x)=r\Big(\frac{x-x_0}{\eps}\Big)
e^{\frac{\iu}{\eps}x\cdot\xi_0},\qquad x\in\R^N.
\end{equation}
We shall assume that $N\geq 1$, $0<p<2/N$,
$\iu$ is the imaginary unit and $r\in H^1\cap C^2(\R^N)$ 
is the unique~\cite{kwong} (up to translations) 
positive and radially symmetric solution of the elliptic problem
\begin{equation}
\label{seMF}
\tag{$E$}
-\frac{1}{2}\Delta r+\lambda r=r^{2p+1}\qquad\text{in $\R^N$},
\end{equation}
for some value $\lambda>0$.\
Finally, $x_0$ and $\xi_0$ are given vectors in $\R^N$ that should be conveniently thought
(in the transition from quantum to classical mechanics) as corresponding to the 
{\em initial position} and {\em initial velocity} respectively of 
a point particle.

In this framework, since~\eqref{probMF} 
has a conservative nature, the typical expected behaviour is 
that the solution travels with the shape of $r((x-x(t))/\eps)$ (hence
its support shrinks, as $\eps$ gets small)
along a suitable concentration line $x(t)$ merely depending on the potential $V$
and starting at $x_0$ with initial slope $\xi_0$.

On the basis of the analytical results 
currently available in literature (see the discussion in Section \ref{theoryfacts}), we believe
that providing some numerical study is useful to complete the overall picture
of this phenomenon and furnish some practical machinery for the computation 
of the solutions of~\eqref{seMF} and, in turn, of~\eqref{probMF}-\eqref{initialD}. 
The authors are not aware of any other contribution
in the literature on this issue. For the linear Schr\"odinger, some results can be found in
\cite{jinyang}.

\subsection{Facts from the theory}
\label{theoryfacts}
It is well-known that, given a positive real number $m$, 
the afore mentioned (ground state) solution $r$ of~\eqref{seMF} (where the value of $\lambda$
depends on $m$) can obtained through
the following variational characterization on the sphere of $L^2(\R^N)$
\begin{equation}
\label{gsminchar}
{\mathcal E}(r)=\inf\{{\mathcal E}(u): u\in H^1(\R^N),\,\,\|u\|_{L^2}^2=m\},
\end{equation}
where ${\mathcal E}:H^1(\R^N)\to\R$ is the $C^2$ energy functional
\begin{equation}
\label{energyfunct}
{\mathcal E}(u)=\frac{1}{2}\int_{\R^N}|\nabla u|^2 dx-\frac{1}{p+1}\int_{\R^N}|u|^{2p+2}dx.
\end{equation}
Furthermore, there exists a suitable choice of $m$
yielding $\lambda=1$ as eigenvalue in equation~\eqref{seMF}.
The restriction to the values of $p$ below $2/N$ is strictly
related to the global well-posedness of~\eqref{probMF} for any choice of initial
data $\phi_0$ in $H^1$. If $p$ is larger than or equal to $2/N$, then the solution can blow-up in
finite time (see e.g.\ the monograph by T.\ Cazenave~\cite{cazenave}). In particular, in the two 
dimensional case, $p$ will be picked in $(0,1)$.

From the analytical side, it has been rigorously known since 2000 that
the solution $\phi_\eps(t)$ of~\eqref{probMF} remains close to the ground state
$r$, in the sense stated here below, locally uniformly in time,
as $\eps$ is pushed to zero. 
As we said, this dynamical behaviour is typically known as soliton dynamics (for a 
recent general survey on solitons and their stability, see the work of T.\ Tao~\cite{tao-solitons}). 

For the nonlinear equation~\eqref{probMF}, rigorous results about 
the soliton dynamics were obtained in various papers
by J.C.\ Bronski, R.L.\ Jerrard~\cite{bronski} and S.\ Keraani~\cite{Keerani1,Keerani2}.
We also refer to~\cite{squamagn} for a complete study of the problem with the additional
presence of an external time independent magnetic vector potential $A:\R^N\to\R^N$,
and to~\cite{mopelsqu} for a study of a system of two coupled nonlinear Schr\"odinger
equations, a topic which is rapidly spreading in the last few years.
The arguments are mainly based upon the following ingredients: 
the energy convexity estimates proved by M.\ Weinstein~\cite{weinsteinMS,weinstein2}
to get the so called modulational stability, the use of conservation 
laws (mass and energy) satisfied by the equation and by the associated
Hamiltonian system in $\R^N$ built upon the guiding external 
potential $V$, that is the classical {\em Newton law}
\begin{equation}
	\label{newtonsenza}
\begin{cases}
\ddot x(t)=-\nabla V(x(t)),  &\\
     x(0)=x_0, &\\
\dot x(0)=\xi_0.
\end{cases}
\end{equation}
Under reasonable assumptions on $V$ (e.g.\ uniform boundedness of the second order partial 
derivatives), equation \eqref{newtonsenza}
admits a unique global solution $(x(t),\xi(t))$ which satisfies
the following conservation law
$$
{\mathcal H}(t)=\frac{1}{2}|\xi(t)|^2+V(x(t)),\qquad {\mathcal H}(t)={\mathcal H}(0),\quad t\geq 0.
$$
Let us now define a suitable scaling of the standard norm of $H^1(\R^N)$
$$
\|\phi\|_{\H_\eps}^2=\eps^{2-N}\|\nabla \phi\|_{L^2}^2+\eps^{-N}\|\phi\|_{L^2}^2,\qquad\eps>0.
$$
\vskip5pt
\noindent
The precise statement of the soliton dynamics reads as follows
\begin{ana}[Cf.\ e.g.\ \cite{bronski,Keerani2}]
\label{anal11}
Let $\phi_\eps(t)$ be the solution to problem \eqref{probMF}
corresponding to the initial datum~\eqref{initialD}.
Then there exists a family of shifts $\teta_{\vep}:\R^+\to [0,2\pi)$ such that,
as $\eps$ goes to zero, $\phi_\eps(x,t)$ is equal to the function 
\begin{equation}
\label{mainconclusintro}
\phi_\eps^r(x,t)=
r\Big(\frac{x-x(t)}{\vep}\Big)e^{\frac{{\rm i}}{\vep}\left[x\cdot  \dot x(t)+\teta_{\vep}(t)\right]},
\qquad x\in\R^N,\,\,t>0,
\end{equation}
up to an error function $\omega_\eps(x,t)$ such that
$\|\omega_\eps(t)\|_{\H_\eps}\leq {\mathcal O}(\eps)$,
locally uniformly in time.
\end{ana}
It is important to stress that, in the particular case of {\em standing wave solutions}  of~\eqref{probMF}, namely special solutions of~\eqref{probMF} of the form 
$$
\phi_\eps(x,t)=u_\eps(x) e^{-\frac{\iu}{\eps}\theta t},\quad x\in\R^N,\, t\in\R^+,\quad(\theta\in\R),
$$ 
where $u_\eps$ is a real-valued function, there is an enormous
literature regarding the semi-classical limit for
the corresponding elliptic equation
$$
-\textstyle{\frac{\eps^2}{2}}\Delta u_\eps+V(x)u_\eps=|u_\eps|^{2p}u_\eps,\quad x\in\R^N.
$$
See the recent book~\cite{ambook} by A.\ Ambrosetti and A.\ Malchiodi and the references therein. 
To this regard notice that, if
$\xi_0=0$ (null initial velocity) and $x_0$ is a critical point of the 
potential $V$, as equation~\eqref{newtonsenza} admits the trivial solution $x(t)=x_0$ and $\dot x(t)=0$
for all $t\in\R^+$, formula~\eqref{mainconclusintro} reduces to
\begin{equation*}
\phi_\eps^r(x,t)=r\Big(\frac{x-x_0}{\vep}\Big)
e^{\frac{{\rm i}}{\vep}\teta_{\vep}(t)},\qquad x\in\R^N,\,\,t>0,
\end{equation*}
so that the concentration of $\phi_\eps(t)$ is {\em static}
and takes place at $x_0$, instead occurring 
along a smooth concentration curve in $\R^N$. 
This is consistent
with the literature for the standing wave solutions mentioned above.

For other achievements about the full dynamics of~\eqref{probMF}, see also~\cite{gril1,gril2} 
(in the framework of orbital stability of standing waves) as well as~\cite{kaup,keener} (in the framework
of non-integrable perturbation of integrable systems). Similar results were investigated in geometric optics 
by a different technique (WKB method), namely
writing formally the solution as $u_\eps=U_\eps(x,t) e^{\iu\theta(x,t)/\eps}$, with
$U_\eps=U_0
+\eps U_1+\eps^2 U_2\cdots,$ 
where $\theta$ and $U_j$ are solutions, respectively, of a Hamilton-Jacobi type
equation (the eikonal equation) and of a system of transport equations.
In presence of a constant external potential, the orbital stability issue for 
problem~\eqref{probMF} was investigated by T.\ Cazenave and P.L.\ Lions~\cite{cl}, and by
M.\ Weinstein in~\cite{weinsteinMS,weinstein2}. Then, A.\ Soffer and M.\ Weinstein proved in~\cite{soffer1}
the asymptotic stability of nonlinear ground states of~\eqref{probMF}. 
See also the following important contributions: 
V.\ Buslaev and G.\ Perelman~\cite{buslaev1},  
V.\ Buslaev and C.\ Sulem~\cite{buslaev3},
J.\ Fr\"ohlich, S.\ Gustafson, L.\ Jonsson, I.M.\ Sigal, T.-P.\ Tsai and H.-T.\ Yau~\cite{frolich1,frohl3,aboufrosig},
J.\ Holmer and Zworski~\cite{holmer},
A.\ Soffer and M.\ Weinstein~\cite{soffer2,soffer3},
T.-P.\ Tsai and H.-T.\ Yau~\cite{tsai1}.

\vskip2pt
\noindent
Another interesting problem concerns the case where the initial datum
is multibump (for simplicity two bumps), say,
\begin{equation}
\label{initialD2}
\phi_0(x)=r_1\Big(\frac{x-x_0}{\eps}\Big)
e^{\frac{\iu}{\eps}x\cdot\xi_0}+
r_2\Big(\frac{x-y_0}{\eps}\Big)
e^{\frac{\iu}{\eps}x\cdot\eta_0},\qquad x\in\R^N.
\end{equation}
where $r_i$ are solutions to the problem
$$
{\mathcal E}(r_i)=\inf\{{\mathcal E}(u): u\in H^1(\R^N),\,\,\|u\|_{L^2}^2=m_i\},
$$
for some fixed $m_i>0$, $i=1,2$ and $x_0,y_0,\xi_0,\eta_0$ are taken as initial data for 
\begin{equation}
	\label{newtonsenza22}
\begin{cases}
\ddot x(t)=-\nabla V(x(t)),  &\\
     x(0)=x_0, &\\
\dot x(0)=\xi_0.
\end{cases}
\qquad
\begin{cases}
\ddot y(t)=-\nabla V(y(t)),  &\\
     y(0)=y_0, &\\
\dot y(0)=\eta_0.
\end{cases}
\end{equation}

Then we state the following

\begin{ana}[Cf.\ e.g.\ \cite{aboufrosig}]
Let $\phi_\eps(t)$ be the solution to~\eqref{probMF}
corresponding to the initial datum~\eqref{initialD2}.
Then there exist two families of shifts $\teta_{\vep}^i:\R^+\to [0,2\pi)$ such that,
as $\eps$ goes to zero, $\phi_\eps(x,t)$ is equal to the function 
\begin{equation}
\label{mainconclusintro2}
\phi_\eps^r(x,t)=
r_1\Big(\frac{x-x(t)}{\vep}\Big)e^{\frac{{\rm i}}{\vep}\left[x\cdot  \dot x(t)+\teta_{\vep}^1(t)\right]}
+
r_2\Big(\frac{x-y(t)}{\vep}\Big)e^{\frac{{\rm i}}{\vep}\left[x\cdot  \dot y(t)+\teta_{\vep}^2(t)\right]},
\end{equation}
up to an error function $\omega_\eps(x,t)$ depending both on $\eps$ and on the initial
relative velocity $v=|\xi_0-\eta_0|$ (the larger is $v$ the smaller is the error), 
locally uniformly in time.
\end{ana}

\noindent
See figure \ref{fig:2} in the final section for a movie showing this behaviour.
\medskip

\section{Numerical computation of the soliton dynamics}
In the numerical simulations included in the last section of the paper, we shall
consider the two dimensional case. On the other hand, here we consider the general case.

\subsection{Overview of the method}
\begin{figure}
\includegraphics[scale=0.6]{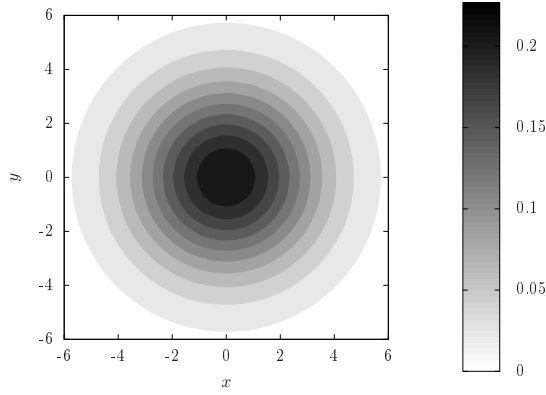}
\caption{The positive, radially symmetric and radially decreasing 
ground state solution $r_\infty$ of~\eqref{gsminchar} with $m=1$ and $p = 0.2$. 
Of course, in the computation of $r_\infty$, there is a spurious imaginary part of maximum value around $10^{-16}$, since the complex FFT algorithm is involved. The corresponding value of
$\lambda_\infty$ is $\lambda_\infty = -0.37921$.}
\label{fig:0}
\end{figure}

Our purpose is to solve the Schr\"odinger equation
\begin{equation}\label{eq:schroedinger}
\left\{\begin{aligned}
&\iu \eps \partial_t \phi_\varepsilon(x,t)=
-\frac{\varepsilon^2}{2}\Delta\phi_\varepsilon(x,t)+
V(x)\phi_\varepsilon(x,t)-\abs{\phi_\varepsilon(x,t)}^{2p}\phi_\varepsilon(x,t),&&
x\in\RN,\\
&\phi_\varepsilon(x,0)=r_\varepsilon(x-x_0),&&x\in\RN,
\end{aligned}\right.
\end{equation}
where $r_\varepsilon(x)=u(x/\varepsilon)$, 
so that $\phi(x,t)=u(x)e^{-\iu \lambda t}$ is the solution of 
\begin{equation}
	\label{eq:gs}
\iu \partial_t \phi(x,t)=-\frac{1}{2}\Delta \phi(x,t)-
\abs{\phi(x,t)}^{2p}\phi(x,t)
\end{equation}
being $u$ real, positive
and minimizing the energy~\eqref{energyfunct}
under the constraint $\norm[L^2]{u}^2=m$.
Instead of a direct minimization of the energy (see, e.g., \cite{BT03,CORT09}), 
here we consider
the following parabolic differential equation
\begin{equation}\label{eq:parabolic}
\left\{
\begin{aligned}
&\partial_t r(x,t)=\frac{1}{2}\Delta r(x,t)+r^{2p+1}(x,t)+
\lambda(r(x,t))r(x,t),&& x\in\RN,\ t>0\\
&r(x,0)=r_0(x),\ \norm[L^2]{r_0}^2=m,&& x\in\RN
\end{aligned}\right.
\end{equation}
with vanishing boundary conditions, where the map $t\mapsto \lambda(r(\cdot,t))$ is defined by
\begin{equation*}
\lambda(r(x,t))=\frac{\frac{1}{2}\int_{\RN}\abs{\nabla r(x,t)}^2\rmd x-
\int_{\RN}\abs{r(x,t)}^{2p+2}\rmd x}{\norm[L^2]{r}^2}
\end{equation*}
This approach is similar to the imaginary time method (see, e.g., \cite{BD04}),
based on the propagation
of the Schr\"odinger equation along imaginary time $-\iu t$ and projection
to the $L^2$ sphere of radius $\sqrt{m}$.
In equation~\eqref{eq:parabolic}, projection is not necessary and
the energy decreases: in fact,
if we multiply equation~\eqref{eq:parabolic} by $r(x,t)$ and integrate over 
$\RN$, we easily get 
\begin{equation*}
\frac{1}{2}\frac{\rmd}{\rmd t}\norm[L^2]{r(\cdot,t)}^2=\int_\RN r(x,t)\partial_t r(x,t)\rmd x = 0
\end{equation*}
and if we multiply equation~\eqref{eq:parabolic} by $\partial_t r(x,t)$ 
and integrate over  $\RN$, we get
\begin{align*}
\frac{1}{2}\frac{\rmd}{\rmd t} {\mathcal E}(r(\cdot,t))&=
-\int_\RN\abs{\partial_t r(x,t)}^2\rmd x+
\lambda(r(x,t))\int_\RN r(x,t)\partial_t r(x,t)\rmd x   \\
&=
-\int_\RN\abs{\partial_t r(x,t)}^2\rmd x\le0
\end{align*}
Hence, the steady-state solution $r_\infty(x)=r(x,t\to\infty)$ of~\eqref{eq:parabolic}
satisfies $\norm[L^2]{r_\infty}^2=m$ and has a minimal energy. 
In fact, notice that by the results of~\cite{kwong}, 
for any $\lambda>0$ there exists a unique (up to translations)
positive and radially symmetric solution $r=r_\lambda$ of~\eqref{seMF}. In turn,
given $\lambda_1,\lambda_2>0$, if $r_1,r_2:\R^N\to\R$ denote, respectively, the positive radial
solutions of the equations
\begin{equation*}
-\frac{1}{2}\Delta r_1+\lambda_1 r_1=r_1^{2p+1},\qquad\,\,\,
-\frac{1}{2}\Delta r_2+\lambda_2 r_2=r_2^{2p+1},
\end{equation*}
then it is readily verified that 
\begin{equation*}
r_2(x)=\mu r_1(\gamma x),\qquad \gamma=\left(\frac{\lambda_2}{\lambda_1}\right)^\frac{1}{2},\quad
\mu=\left(\frac{\lambda_2}{\lambda_1}\right)^\frac{1}{2p},
\end{equation*}
which tells us that that, up to a scaling, the solution corresponding
to different values of $\lambda$ is unique. 
Notice now that, due to the 
choice of the bump like initial datum 
(Gaussian like, see~\eqref{initialparabolic}) 
in the iterations to compute $r_\infty$ (see the discussion below), it turns out that 
$\lambda_\infty$, defined as $\lambda(r_\infty)$,
is negative and $r_\infty$ is positive, radially symmetric (see figure~\ref{fig:0}) and solves
\begin{equation*}
-\frac{1}{2}\Delta r_\infty+\hat\lambda_\infty r_\infty=r_\infty^{2p+1},
\end{equation*}
where $\hat\lambda_\infty=-\lambda_\infty>0$. If $r_m$ denotes the ground state
solution (with the corresponding positive eigenvalue denoted by $\lambda_m$), then we have
\begin{equation}
\label{realvscomp}
r_m(x)=\mu r_\infty(\gamma x),\qquad \gamma=\left(\frac{\lambda_m}{\lambda_\infty}\right)^\frac{1}{2},\quad
\mu=\left(\frac{\lambda_m}{\lambda_\infty}\right)^\frac{1}{2p}. 
\end{equation}
On the other hand, by construction, we have
$$
m=\|r_m\|_{L^2}^2=\int_\RN r_m^2(x)dx=\mu^2\gamma^{-N}\int_\RN r_\infty^2(x)dx=m\mu^2\gamma^{-N},
$$
namely $\mu^2\gamma^{-N}=1$. Finally, by the definition of $\gamma$ 
and $\mu$ in~\eqref{realvscomp}, we get $\lambda_m=\hat \lambda_\infty$ and $\gamma=\mu=1$,
yielding from~\eqref{realvscomp} the desired conclusion, that is
$$
r_\infty=r_m.
$$
Moreover, $r_\infty(x) e^{-\iu \lambda(r_\infty(x))t}$ is a solution of
\eqref{eq:gs}. 
We will take $r_\varepsilon(x-x_0)=r_\infty((x-x_0)/\varepsilon)$
as our candidate initial condition for the time-dependent nonlinear 
Schr\"odinger equation~\eqref{eq:schroedinger}.
From a numerical point of view,
it is convenient to compute directly $r_\infty(x/\varepsilon)$ instead
of $r_\infty(x)$ and to apply the change of variable 
$\Phi(X,t)=\sqrt[4]{\varepsilon^N}\phi_\varepsilon(x,t)$,
 $\sqrt{\varepsilon}X=x$,
 to the nonlinear Schr\"odinger
equation~\eqref{eq:schroedinger}, and hence to equation~\eqref{eq:parabolic}. 
Altogether, we need to solve
\begin{equation}\label{eq:parabolictosolve}
\left\{
\begin{aligned}
&\partial_t R(X,t)=\frac{\varepsilon}{2}\Delta R(X,t)+
\varepsilon^{-Np/2}R(X,t)^{2p+1}+
\Lambda(R(X,t))R(X,t),&&X\in\RN\\
&R(X,0)=R_0(X),\ \norm[L^2]{R_0}^2=m\varepsilon^N,&&X\in\RN
\end{aligned}\right.
\end{equation}
with
\begin{equation*}
\Lambda(R(X,t))=
\frac{\frac{\varepsilon}{2}\int_{\RN}\abs{\nabla R(X,t)}^2\rmd X-
\varepsilon^{-Np/2}\int_{\RN}\abs{R(X,t)}^{2p+2}\rmd X}{\norm[L^2]{R}^2}
\end{equation*}
where $R(X,t)=r(x/\varepsilon,t)\sqrt[4]{\varepsilon^N}$. Since
it is not possible to numerically integrate the equation up to an
infinite time, we will consider $R(X,\tss)$ the steady-state 
as soon as ${\mathcal E}(R(X,\tss))$ is stabilized within a prescribed tolerance.
The initial
condition $R_0(X)$ can be arbitrarily chosen (in the class of bump like functions), but an initial solution
with small energy will shorten the ``steady-state'' time $\tss$.
Among the family of the 
Gaussian functions parameterized by $\sigma$
\begin{equation}
\label{initialparabolic}
R^{\sigma}(X)=\sqrt{m}\sigma^{N/2}e^{-\abs{\sigma X/\sqrt{\varepsilon}}^2/2}
\sqrt[4]{\left(\frac{\varepsilon}{\pi}\right)^N}
\end{equation}
with $\norm[L^2]{R^\sigma}^2=m\varepsilon^N$ it is possible to choose
the one with minimal energy. In fact
\begin{equation*}
{\mathcal E}(R^\sigma)=\sigma^2\frac{\varepsilon}{2}\int_\RN\abs{\nabla R^1(X)}^2\rmd X
-\sigma^{Np}\frac{\varepsilon^{-Np/2}}{p+1}\int_\RN \abs{R^1(X)}^{2p+2}\rmd X
\end{equation*}
If we define
\begin{equation*}
A=\frac{\varepsilon}{2}\int_\RN \abs{\nabla R^1(X)}^2\rmd X,\quad
B=\frac{\varepsilon^{-Np/2}}{p+1}\int_\RN \abs{R^1(X)}^{2p+2}\rmd X
\end{equation*}
the minimum for ${\mathcal E}(R^\sigma)$ is attained for
\begin{equation*}
\sigma=\left(\frac{BNp}{2A}\right)^{\frac{1}{2-Np}}
\end{equation*}
The quantities $A$ and $B$ can be analytically computed and give
\begin{equation*}
A=\left\{\begin{aligned}
&\frac{m\varepsilon}{4}&&N=1\\
&\frac{m\varepsilon^2}{2}&&N=2\\
&\frac{3m\varepsilon^3}{4}&&N=3
\end{aligned}\right.,\quad B=\frac{m^{p+1}\varepsilon^N}{\pi^{Np/2}(p+1)^{1+N/2}}
\end{equation*}

\subsection{Numerical discretization}
With the normalization introduced above, the nonlinear Schr\"odinger equation 
to solve is
\begin{equation} \label{eq:schroedingertosolve}
\left\{\begin{aligned}
&\iu\partial_t\Phi(X,t)=
-\frac{1}{2}\Delta\Phi(X,t)+
\frac{V(\sqrt{\varepsilon}X)}{\varepsilon}\Phi(X,t)
-\frac{\abs{\Phi(X,t)}^{2p}}{\sqrt{\varepsilon^{2+Np}}}\Phi(X,t),&&
X\in\RN\\
&\Phi(X,0)=R(X-X_0,\tss),&&X\in\RN
\end{aligned}\right.
\end{equation}
A well-established numerical method for the cubic Schr\"odinger equation
(focusing or defocusing case) is the Strang splitting \cite{BJM03,BJaM03,CNT09}.
It is based on a split of the full equation into two parts,
in which the first is spectrally discretized in space and then exactly
solved in time and the second has an analytical solution. 
We used the Strang splitting method as well.
The first part is
\begin{subequations}
\begin{equation}\label{eq:schro1}
\iu\partial_t\Phi_1(X,t)=
-\frac{1}{2}\Delta\Phi_1(X,t)
\end{equation}
Thus, the Fourier coefficients of $\Phi_1(X,t)$ restricted to a sufficiently
large space domain satisfy a linear and 
diagonal system of ODEs, which can be exactly solved.
The second part is
\begin{equation}\label{eq:schro2}
\iu\partial_t\Phi_2(X,t)=
\frac{V(\sqrt{\varepsilon}X)}{\varepsilon}\Phi_2(X,t)
-\frac{\abs{\Phi_2(X,t)}^{2p}}{\sqrt{\varepsilon^{2+Np}}}\Phi_2(X,t)
\end{equation}
\end{subequations}
It is easy to show that the quantity $\abs{\Phi_2(X,t)}^{2p}$ is constant in 
time for this equation. Then it has an analytical solution.
Given the approximated solution $\Phi_n(X)\approx\Phi(X,t_n)$ of 
equation~\eqref{eq:schroedingertosolve}, a single time step of 
the Strang splitting Fourier spectral method can be summarized in
\begin{enumerate}
\item take $\Phi_n(X)$ as initial
solution at time $t_n$ for~\eqref{eq:schro1} and solve 
for a time step $k/2$, obtaining 
$\Phi_1(X,t_n+k/2)$;
\item take $\Phi_1(X,t_n+k/2)$ as initial solution at time $t_n$ for 
\eqref{eq:schro2} and solve for a time step $k$, obtaining
$\Phi_2(X,t_n+k)$;\label{realspace}
\item take $\Phi_2(X,t_n+k)$
as initial solution for~\eqref{eq:schro1} and solve for a time step $k/2$,
obtaining $\Phi_{n+1}(X)$.
\end{enumerate}
The result $\Phi_{n+1}(X)$ is an approximation of $\Phi(X,t_n+k)$. 
Since the solutions of the first part and the second part are trivial to
compute in
the spectral space and in the real space, respectively, it is necessary
to transform the solution from spectral space to real and
from real space to spectral before and 
after step~\eqref{realspace} above, respectively.
All the transformations can be carried out by the FFT algorithm.
The method turns out
to be spectrally accurate in space and of the second order in time.

Therefore, we used the Fourier spectral decomposition 
for the solution
of equation~\eqref{eq:parabolictosolve}, too. Together with the Galerkin
method, it yields a nonlinear system of ODEs
\begin{equation}\label{eq:ODEs}
\left\{\begin{aligned}
&\hat R'(t)=\frac{\varepsilon}{2}D \hat R(t)+f(\hat R(t)),&&t>0,\\
&\hat R(0)=\hat R_0,
\end{aligned}\right.
\end{equation}
where $\hat R$ is the vector of Fourier coefficients, $D$ the diagonal
matrix of the eigenvalues of the Laplace operator and $f$ the truncated
Fourier expansion of the whole nonlinear part of
equation~\eqref{eq:parabolictosolve}.
For the solution of equation~\eqref{eq:ODEs} we used an exponential Runge--Kutta
method of order two (see, e.g., \cite{HO05}), 
with the embedded exponential Euler method.
Given the approximation $\hat R_n\approx \hat R(t_n)$, a single time
step of the method is
\begin{enumerate}
\item set $A_{n+1}=k_{n+1}\frac{\varepsilon}{2}D$ and $R_{n1}=R_n$;
\item compute $\hat R_{n2}=\exp(A_{n+1})\hat R_{n1}+
k_{n+1}\varphi_1(A_{n+1})f(\hat R_{n1})$ (exponential Euler method);
\item compute $\hat R_{n+1}=\hat R_{n2}+k_{n+1}\varphi_2(A_{n+1})(-f(\hat R_{n1})+
f(\hat R_{n2}))$\label{error}
\end{enumerate}
where $\varphi_1(z)$ and $\varphi_2(z)$ are the analytic functions
\begin{align*}
\varphi_1(z)&=\frac{e^z-1}{z},\ z\ne 0, &  
\varphi_2(z)&=\frac{e^z-1-z}{z^2},\ z \ne 0,\\
\varphi_1(0)&=1,               & \varphi_2(0)&=\frac{1}{2}.
\end{align*}
The result is an approximation of $\hat R(t_n+k_{n+1})$.
Exponential integrators are explicit and do not suffer of time step
restrictions. However, they require the computation of matrix functions. 
In our case, the matrices involved $A_{n+1}$ are diagonal and the computation
of the matrix functions $\exp(A_{n+1})$, $\varphi_1(A_{n+1})$ and 
$\varphi_2(A_{n+1})$ is trivial.
In order to compute the terms $f(\hat R_{n1})$ and $f(\hat R_{n2})$, it 
is necessary to recover the functions in the real space corresponding 
to the Fourier spectral coefficients $\hat R_{n1}$ and $\hat R_{n2}$,
respectively, then compute the nonlinear part of 
equation~\eqref{eq:parabolictosolve} and finally to compute its Fourier 
transform. All the transformations can be carried out by the FFT algorithm.
The term $\hat R_{n+1}-\hat R_{n2}$ in step~(\ref{error}) above
can be used as an error estimate for $R(t_{n+1})-R_{n+1}$ 
and then it is possible
to derive a variable time step integrator. This is particularly
useful for our aim of computing the steady-state of the equation: in fact,
we expect that the as soon as the solution approaches the steady-state it is
possible to enlarge the time step, thus reducing the computational cost.
The method turns out to be spectrally accurate in space and of the second
order in time.

\bigskip
\section{Two dimensional examples and error analysis}

\begin{figure}
\centering
\href{run:movie1bh.avi}{\includegraphics[scale=0.8]{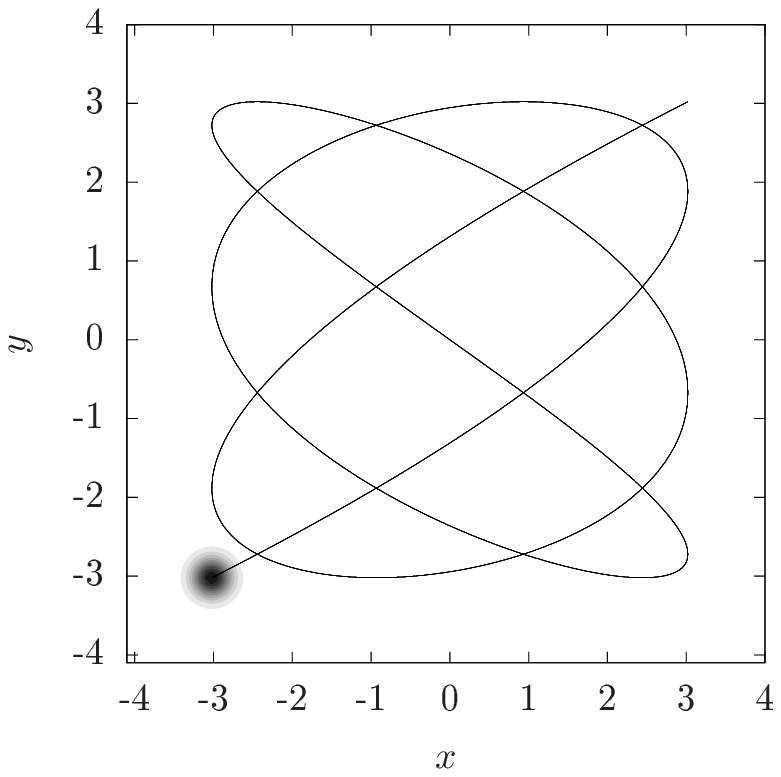}}
\hspace{1cm}
\href{run:movie1bhirr.avi}{\includegraphics[scale=0.8]{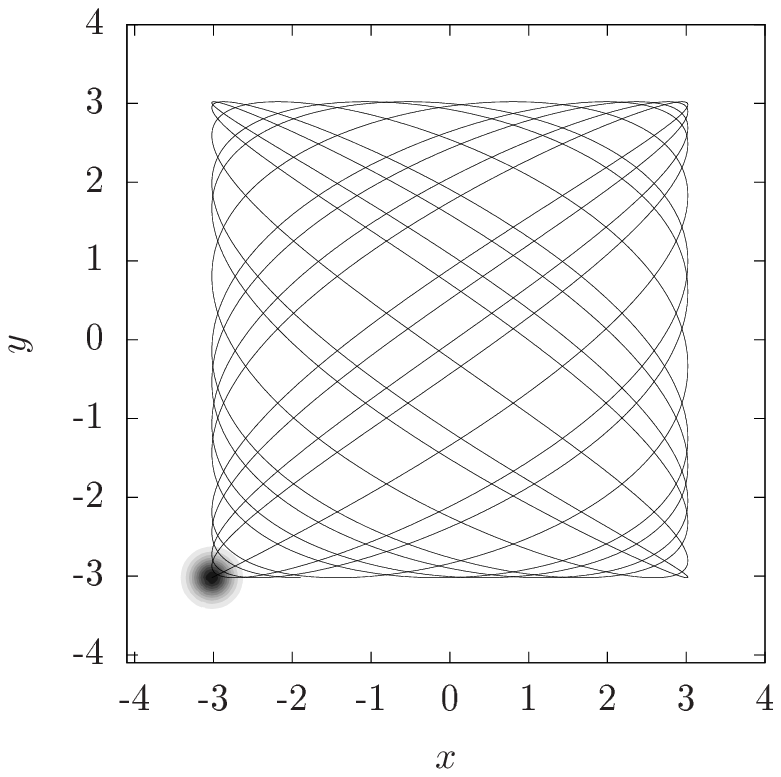}}
\caption{In both simulation movies we set $\eps=0.01$, $p = 0.02$, $m = 1$,
$(x_0,y_0)=(-3.0,-3.0)$, $v_0 = (0,0)$. In the left movie, we chose
$\omega_1=1.4$ and $\omega_2=1$ (rational ratio). In the right movie,
we chose $\omega_1=\sqrt{2}$ and $\omega_2=1$ (irrational ratio). Notice that,
although the ratios $\omega_2/\omega_1$ are very close in the two examples,
the soliton dynamics is ergodic in the right movie. Of course 
the figures refer to the (squared modulus of the) solution
at the time $t=0$ and contain the concentration paths 
(admitting an analytic expression) that the soliton is going
to travel on.}
\label{fig:1}
\end{figure}

\begin{figure}
\centering
\href{run:movie2bh.avi}{\includegraphics[scale=0.8]{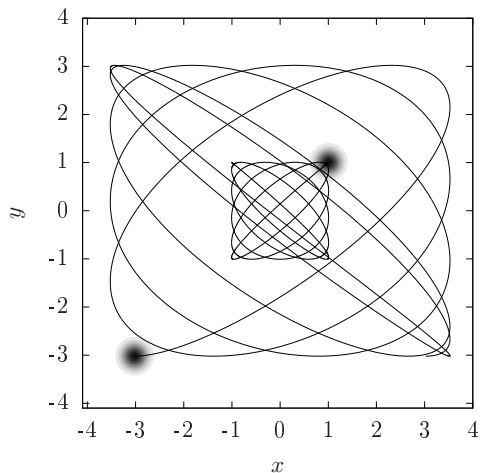}}
\caption{In the simulation movie, we set $\eps=0.01$, $p = 0.02$, $m = 1$,
$(x_0^1,y_0^1)=(-3,-3)$, $(x_0^2,y_0^2)=(1,1)$, $v_0^1=(2,0)$,
$v_0^2=(0,0)$, $\omega_1=1.1$ and $\omega_2=1$.}
\label{fig:2}
\end{figure}

\begin{figure}
\centering
\includegraphics[scale=0.75]{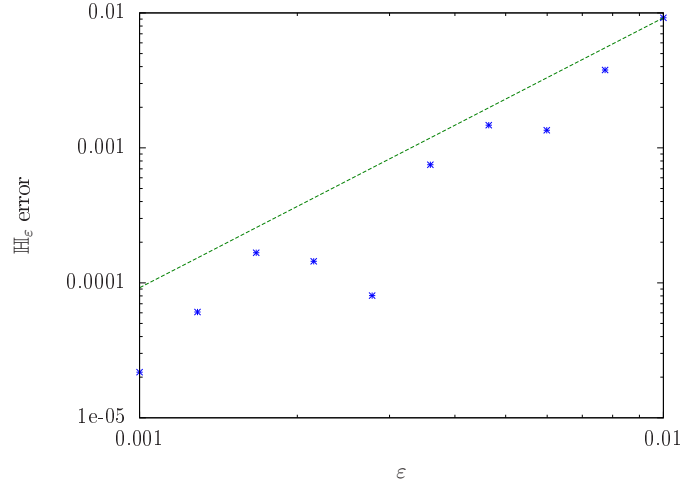}
\caption{For the error analysis, we set $p=0.02$, $(x_0,y_0)=(-0.5,-0.5)$, $m=1$, $\omega=(2,1)$
and a final time $t=\pi$.
With the change of variable we used, $\sqrt{\varepsilon}X=x$ and 
$U(X)=\sqrt{\varepsilon}u(x)$, we have
$\lVert u\rVert_{L^2}^2=\lVert U\rVert_{L^2}^2$ and
$\lVert \nabla_x u\rVert_{L^2}^2=\frac{1}{\varepsilon}\lVert \nabla_X U\rVert_{L^2}^2$.
Hence, the numerical error is computed through formula (written for the 2D case)
$\| u\|_{\H_\varepsilon}=
\sqrt{\varepsilon^{-1}\| \nabla_X U\|_{L^2}^2+
\varepsilon^{-2}\| U\|_{L^2}^2}$.
As predicted by the analytical property~\ref{anal11}, the error in the 
$\| \cdot\|_{\H_\varepsilon}$ is below the order ${\mathcal O}(\eps)$.}
\label{fig:3}
\end{figure}

In this section, in order to provide some examples, 
we reduce to the two dimensional setting and focus on 
the physically relevant case of harmonic potential
\begin{equation*}
V(x,y)=\omega_1^2 x^2+\omega_2^2 y^2,\qquad
\omega_1,\omega_2>0,
\end{equation*}
well-established in the theory of Bose-Einstein condensates.
In the two movies starting from figure~\ref{fig:1} we show the dynamics of the solitary wave
along two Lissajous curves, periodic in the left side and ergodic for the right side.
In the movie starting from figure~\ref{fig:2} we report the soliton dynamics in the
case of an initial datum exhibiting a double bump behaviour (with a sufficiently
large distance between the centers) up to the collision time.
It is important to stress that in these figures the paths have an
analytic expression and are plotted before the dynamics starts. The movies will then
show that the centers of mass of the solitons follow adherently these curves up to the final
computation time. An analysis of the error (in the single bump case) arising when the modulus of the solution
$|\phi_\eps(x,t)|$ is replaced by the modulus of the expression in the representation 
formula~\eqref{mainconclusintro}, namely $r((x-x(t))/\eps)$, is indicated in figure~\ref{fig:3}.
As predicted by the analytical property~\ref{anal11}, the error in the 
$\| \cdot\|_{\H_\varepsilon}$ is below the order ${\mathcal O}(\eps)$.

\bigskip
\bigskip

\end{document}